\documentclass[11pt,a4paper]{amsart}
\usepackage{amssymb,amsmath,a4}
\usepackage[T1]{fontenc} 
\allowdisplaybreaks
\def\R#1#2#3#4{\hat{R}^{#1#2}_{#3#4}}
\def\Rm#1#2#3#4{\hat{R}^-{}^{#1#2}_{#3#4}}
\def\Rl#1#2#3#4{\grave{R}^{#1#2}_{#3#4}}
\def\Rlm#1#2#3#4{\grave{R}^-{}^{#1#2}_{#3#4}}
\def\Rr#1#2#3#4{\acute{R}^{#1#2}_{#3#4}}
\def\Rrm#1#2#3#4{\acute{R}^-{}^{#1#2}_{#3#4}}
\def\Rc#1#2#3#4{\check{R}^{#1#2}_{#3#4}}
\def\Rcm#1#2#3#4{\check{R}^-{}^{#1#2}_{#3#4}}
\def\d{\mathrm{d}}
\def\id{\mathrm{id}}

\def\real{\mathrm{I\kern-.24em R}}
\def\SUq{\mathrm{SU}_q}
\def\Sq{\mathrm{S}_q}
\def\DeltaR{\Delta_\mathrm{R}}
\def\PhiR{\Phi_\mathrm{R}}
\renewcommand{\H}[2][\mathrm{S}]{\Omega^{#1}_{#2}}
\def\sumL#1{\sum\limits_{#1=1}^N q^{-2#1}}
\def\sumR#1{\sum\limits_{#1=1}^N}
\def\smallsumL#1{{\textstyle\sum\limits_{#1} q^{-2#1}}}
\def\smallsumR#1{{\textstyle\sum\limits_{#1}}}
\let\nob=\nobreakdash
\theoremstyle{plain}
\newtheorem{Prop}{Proposition}
\newtheorem{Thm}[Prop]{Theorem}

\theoremstyle{definition}

\newcommand{\cA}{\mathcal{A}}
\def\nat{\mathrm{I\kern-.18em N}}
\newcommand{\CPq}[2][]{{}^{#1}\mathrm{CP}_q^{#2}}
\newcommand{\x}[2]{x_{#1#2}}
\newcommand{\HCP}{{\,\underline{\smash{\!\H[]{}\!}}\,}}
\newcommand{\RCP}[8]{{\,\underline{\!\hat{R}\!}\,}^{#1#2#3#4}_{#5#6#7#8}}
\newcommand{\RCPc}[8]{{\,\underline{\!\check{R}\!}\,}^{#1#2#3#4}_{#5#6#7#8}}
\newcommand{\RCPm}[8]{{\,\underline{\!\hat{R}\!}\,}^-{}^{#1#2#3#4}_{#5#6#7#8}}
\newcommand{\RCPcm}[8]{%
  {\,\underline{\!\check{R}\!}\,}^-{}^{#1#2#3#4}_{#5#6#7#8}}
\newcommand{\GCP}{{\,\underline{\!\Gamma\!}\,}}
\newcommand{\GCPd}{(\GCP,\d)}
\newcommand{\GCPt}{{\,\underline{\!\tilde{\Gamma}\!}\,}}
\newcommand{\GCPtd}{(\GCPt,\d)}
\newcommand{\GCPtt}{{\,\underline{\!\Tilde{\Tilde{\Gamma}}\!}\,}}
\newcommand{\GCPttd}{(\GCPtt,\d)}
\def\s#1{\mathfrak{s}_{#1}}
\def\si#1{{\mathfrak{s}^{\prime}_{#1}}}
\def\sii#1{{\mathfrak{s}^{\prime\prime}_{#1}}}
\def\siii#1{{\mathfrak{s}^{\prime\prime\prime}_{#1}}}
\def\siv#1{{\mathfrak{s}^{\mathsf{IV}}_{#1}}}
\newcommand{\Gammad}{{(\Gamma,\d)}}

\newcommand{\Gammaat}[1][\alpha\tau]{\Gamma_{#1}}
\newcommand{\Gammaap}[1][\alpha\omega]{\Gamma'_{#1}}
\newcommand{\Gammap}[1][\omega\psi]{\Gamma^{\prime\prime}_{#1}}
\newcommand{\Gammapp}[1][\varrho\tau]{\Gamma^{\prime\prime\prime}_{#1}}
\newcommand{\Gammatl}[1][\lambda]{\tilde{\Gamma}_{#1}}
\newcommand{\Gammatp}[1][\lambda]{\tilde{\Gamma}^{\prime}_{#1}}
\newcommand{\Gammatpp}[1][\lambda]{\tilde{\Gamma}^{\prime\prime}_{#1}}
\newcommand{\Gammato}[1][\lambda]{\tilde{\Gamma}^{\bullet}_{#1}}
\newcommand{\Gammatoo}[1][\lambda]{\tilde{\Gamma}^{\bullet\bullet}_{#1}}
\def\ylth{0.5pt}\def\yubh{1.5ex}\newdimen\yu\yu1ex\newcommand{\yttv}[9]{%
\raisebox{\yubh}{\smash{\rule[0\yu]{#1\yu}{\ylth}\kern-#1\yu\rule[-1\yu]%
{#1\yu}{\ylth}\kern-#1\yu\rule[-2\yu]{#2\yu}{\ylth}\kern-#2\yu\rule[-3\yu]%
{#3\yu}{\ylth}\kern-#3\yu\rule[-4\yu]{#4\yu}{\ylth}\kern-#4\yu}\rule[-#5\yu]%
{\ylth}{#5\yu}\kern-\ylth\kern1\yu\rule[-#5\yu]{\ylth}{#5\yu}\kern-\ylth%
\kern1\yu\rule[-#6\yu]{\ylth}{#6\yu}\kern-\ylth\kern1\yu\rule[-#7\yu]{\ylth}%
{#7\yu}\kern-\ylth\kern1\yu\rule[-#8\yu]{\ylth}{#8\yu}\kern-\ylth\kern1\yu%
\rule[-#9\yu]{\ylth}{#9\yu}\kern-\ylth\kern-5\yu\kern#1\yu\rule[0\yu]{\ylth}%
{\ylth}\kern-\ylth}}\newcommand{\ytt}[9]{~\yttv#1#2#3#4#5#6#7#8#9~}
\newcommand{\ynull}{\textbf{(0)}}
\begin{document}
\author{Martin Welk}
\address{Martin Welk, Department of Mathematics,
University of Leipzig, Augustusplatz~10, 04109 Leipzig,
Germany}
\title[Differential Calculus on Quantum Projective Spaces]%
{Covariant First Order Differential Calculus on Quantum Projective Spaces}
\begin{abstract}
We investigate covariant first order differential calculi on
the quantum complex projective spaces $\CPq{N-1}$ which
are quantum homogeneous spaces for the quantum group $\SUq(N)$.
Hereby, one more well-studied example of covariant first order differential
calculus on a quantum homogeneous space is given.
Since the complex projective spaces are subalgebras of the quantum
spheres $\Sq^{2N-1}$ introduced by Vaksman and Soibelman, we get also
an example of the relations between covariant differential calculus
on two closely related quantum spaces.

Two approaches are combined in obtaining covariant first
order differential calculi on $\CPq{N-1}$: 1.\ restriction of
covariant first order differential calculi from $\Sq^{2N-1}$;
2.\ classification of calculi under appropriate constraints,
using methods from representation theory.

The main result is that under three reasonable settings of dimension
constraints, covariant first order differential calculi on $\CPq{N-1}$
exist and are (for $N\ge6$) uniquely determined.
This is a clear difference as compared to the case of the quantum spheres
where several parametrical series of calculi exist.
For two of the constraint settings, the covariant first order
calculi on $\CPq{N-1}$ are also obtained by restriction from calculi on
$\Sq^{2N-1}$ as well as from calculi on the quantum group $\SUq(N)$.
\end{abstract}
\maketitle
\section{Introduction}
During the last decade, covariant differential calculus on quantum groups has
been under intensive investigation. One main reason for this interest is
that quantum groups are important examples of noncommutative geometric
spaces, equipped with a rich additional structure. The description of
differential calculus on them is an indispensable prerequisite for any
analysis of their geometric structure.

Fundamental concepts of covariant differential calculus on quantum
groups have been introduced in the work of Woronowicz \cite{wr}.
Covariant first order differential calculi on quantum groups have been
constructed; bicovariant first order differential calculi on the most
important quantum groups have been classified \cite{ss1}, \cite{ss2}.
Higher order differential calculus has been studied \cite{s3}, and basic
concepts of differential geometry on quantum groups have been established.

Quantum spaces for quantum groups, and in particular quantum homogeneous
spaces, are a wider class of noncommutative geometric spaces which still
have a rich algebraic structure that can hoped to be helpful in investigating
their geometric properties; they are in many aspects close to quantum groups.
However, unlike for quantum groups, one is still far from having a
comprehensive view on differential calculus on quantum homogeneous spaces.
Only a small number of examples have been studied in detail---apart from
quantum vector spaces which have a comparatively simple structure, we
should mention firstly Podle{\'s}' quantum spheres for which a classification
of covariant first order differential calculi was given in \cite{ap}.
First order differential calculi on the quantum
spheres $\Sq^{2N-1}$ as introduced by Vaksman and Soibelman \cite{vs}
has been classified by the author in a previous paper \cite{wk1}, \cite{wk2}.
The construction of these calculi by restricting covariant calculi from
the quantum group $\SUq(N)$ has been investigated in~\cite{sch}.

In this paper, we want to study the quantum projective spaces $\CPq{N-1}$
which are, as $\Sq^{2N-1}$, quantum homogeneous spaces for the quantum group
$\SUq(N)$, and which are in close relation to the quantum spheres.
With this investigation, one further well-understood example of covariant
first order differential calculus on a quantum homogeneous space is provided;
moreover, the close relationship between $\CPq{N-1}$ and $\Sq^{2N-1}$
allows to ask how this relationship is reflected in the differential
calculus. It is hoped that this work will contribute to a deeper understanding
of covariant differential calculus on quantum homogeneous spaces since
the study of a number of particular examples with
different properties forms the ground on which more general
constructions for differential calculi on quantum spaces and further
theoretical work could be based.

\section{The Quantum Projective Spaces $\CPq{N-1}$}
To start with, we need to recall some basic definitions on quantum spaces.
Our terminology essentially follows~\cite{ks}, \cite{ap}.

Suppose $\cA$ is a Hopf algebra with comultiplication $\Delta$ and counit
$\varepsilon$. A pair $(X,\DeltaR)$ consisting of a unital algebra $X$ and
an algebra homomorphism $\DeltaR: X\to X\otimes\cA$ is called a
\emph{quantum space for $\cA$} if $(\DeltaR\otimes\id)\DeltaR
=(\id\otimes\Delta)\DeltaR$ and $(\id\otimes\varepsilon)\DeltaR
=\id$. Then, $\DeltaR$ is called \emph{(right) coaction} of $\cA$ on $X$.
A quantum space $(X,\DeltaR)$ (or simply, $X$) for $\cA$ is a
\emph{quantum homogeneous space for $\cA$} if there exists an embedding
$\iota:X\to\cA$ with $\DeltaR=\Delta\circ\iota$, i.\,e.\ $X$ can be
considered as a sub-algebra of $\cA$ with the coaction being the restricted
comultiplication.

For our example, the deformation parameter $q$ will always be a real number,
$q\ne 0,\pm1$. Further, $N$ is a natural number, $N\ge2$, parametrising
the dimension of the underlying quantum group $\SUq(N)$.
Throughout the following, the R\nob-matrix $\R{}{}{}{}$ which describes the
commutation relations of the quantum group $\SUq(N)$ \cite{ks}, \cite{rtf},
\cite{ss1} will play an important role. This is an invertible
${N^2\times N^2}$-matrix with $\Rm{}{}{}{}:=\R{}{}{}{}-(q-q^{-1})I$ as its
inverse, where $I$ is the ${N^2\times N^2}$ unit matrix.
The entries of $\R{}{}{}{}$ are given by
\begin{equation*}
\R ijkl = \begin{cases}
1&\text{for $i=l\ne k=j$,}\\
q&\text{for $i=j=k=l$,}\\
q-q^{-1}&\text{for $i=k<j=l$,}\\
0&\text{otherwise.}
\end{cases}
\end{equation*}
For abbreviation, we shall also use the following matrices which are derived
from these fundamental ones:
\begin{align*}
\Rc ijkl&:= \R lkji;&
\Rl ijkl&:= q^{2l-2i}\R jlik;&
\Rr ijkl&:= \R kilj;\\
\Rcm ijkl&:= \Rm lkji;&
\Rlm ijkl&:= q^{2l-2i}\Rm jlik;&
\Rrm ijkl&:= \Rm kilj;
\end{align*}
\begin{align*}
\RCP   stuvijkl &:=\Rlm tuab \R   saic \Rr cbjk \delta_{vl};&
\RCPm  stuvijkl &:=\Rlm tuab \Rm  saic \Rr cbjk \delta_{vl};\\
\RCPc  stuvijkl &:=\Rlm tuab \Rc  bvcl \Rr acjk \delta_{si};&
\RCPcm stuvijkl &:=\Rlm tuab \Rcm bvcl \Rr acjk \delta_{si}.
\end{align*}
Finally, we define 
\begin{align*}
\s+&:=\sum\limits_{i=0}^{N-1}q^{2i};&\si+&:=\s+-1;&
\sii+&:=\si+-q^2;&\siii+&:=\sii+-q^4;&\siv+&:=\siii+-q^6.
\end{align*}

We can now define the quantum projective spaces to be considered here.
These are also given in \cite[11.6]{ks}.
Let $\CPq{N-1}$ be the algebra with $N^2$ generators $\x ij$, $1\le i,j\le N$,
and relations
\begin{align*}
\RCPm stuvijkl \x st\x uv &=q^{-1}\x ij\x kl,&
\RCPc stuvijkl \x st\x uv &=q\x ij\x kl,&
\sumR{i}\x ii&=1.
\end{align*}
The further relation $\sumL{j}\x ij\x jk=\x ik$ is implied by these ones.
The algebra $\CPq{N-1}$ can be equipped with a $*$-structure by letting
$(\x ij)^*:=\x ji$. By the embedding $\iota:\CPq{N-1}\to\SUq(N)$,
$\x ij\mapsto u^1_i(u^1_j)^*=u^1_iS(u^j_1)$ where $u^i_j$ are the $N^2$
coordinates, and $S$ the antipode map of $\SUq(N)$, $\CPq{N-1}$ becomes a
quantum homogeneous space for $\SUq(N)$, with the coaction
$\DeltaR(\x ij)=\sum\limits_{k,l}\x kl\otimes u^k_i S(u^j_l)$.
We shall call $\CPq{N-1}$ \emph{quantum projective space}.

The quantum space $\CPq{N-1}$ can be embedded as a quantum homogeneous
space---i.e.\ respecting its algebra structure and $\SUq(N)$-coalgebra
structure---into the quantum sphere $\Sq^{2N-1}$ which has been
studied in \cite{rtf}, \cite{vs}, \cite{ks}, \cite{wk1}.
The embedding is given by $\x ij\mapsto z_iz^*_j$ where $z_i$, $z^*_i$,
$1\le i\le N$, are the algebra generators of $\Sq^{2N-1}$.

Note that in the literature there exist two different versions of
quantum projective space. The one used here belongs to 
a one-parameter series of quantum projective spaces $\mathcal{B}^\sigma_q$
described in \cite[11.6]{ks}.

\section{First order differential calculus}

The following definitions concerning differential calculus are again
in concordance with, \cite{ks}, \cite{ap}.
By a first order differential calculus on an algebra $X$ we shall mean
a pair $\Gammad$ of a bimodule $\Gamma$ over $X$ and a linear mapping
$\d:X\to\Gamma$ fulfilling Leibniz' rule  $\d(xy)=(\d x)y+x(\d y)$ for all
$x,y\in X$, and $\Gamma=\mathrm{Lin}\{x\d y~ |~ x,y\in X\}$.
The elements of $\Gamma$ are called one-forms.

We call a first order differential calculus $\Gammad$ on a quantum space
$(X,\DeltaR)$ for ${\mathcal A}$ (right) covariant if there exists a
linear mapping $\PhiR:\Gamma\to\Gamma\otimes{\mathcal A}$ which satisfies
the identities $(\PhiR\otimes\id)\PhiR=(\id\otimes\Delta)\PhiR$;
$(\id\otimes\varepsilon) \PhiR=\id$;
$\PhiR(x\omega y)=\DeltaR(x)\PhiR(\omega)\DeltaR(y)$; $\PhiR(\d x)=(\d\otimes
\id)\DeltaR(x)$ for all $x,y\in X$, $\omega\in\Gamma$.
Those one-forms $\omega$ for which the identity
${\DeltaR(\omega)=\omega\otimes1}$ holds are called invariant.
Remember that in the description of bicovariant differential calculi
on quantum groups, invariant one-forms play a central role since essentially
all one-forms can be described by using only left- (or right-) invariant
forms. On quantum homogeneous spaces there exist usually not enough 
invariant one-forms to enable a similar description; nevertheless, they
still are of great significance for the differential calculus.

We mention that if $X$ is a $*$\nob-algebra then the notion of a
$*$\nob-calculus can be introduced: $\Gammad$ is a $*$\nob-calculus
if $\sum\limits_k\,x_k\d y_k=0$ for $x_k,y_k\in X$ always implies
$\sum\limits_k\,\d(y_k^*)x_k^*=0$. However, in this paper $*$\nob-calculi
won't play an important role.

From the definition of a covariant differential calculus and from the
relations of the algebra $\CPq{N-1}$ it is clear that in any
covariant first order differential calculus on $\CPq{N-1}$, the one-form
$\HCP=\displaystyle\sum\limits_{i,j=1}^{N}q^{-2j}\x ij\d\x ji$
is invariant (it might, however, be zero).

\section{Covariant differential calculi on $\CPq{N-1}$}
In this paper, two approaches will be used to obtain covariant first order
differential calculi on the quantum projective spaces.

The first way is based on restricting differential calculi on $\Sq^{2N-1}$
to the sub-algebra $\CPq{N-1}$. Differential calculi on $\Sq^{2N-1}$ have been
classified in~\cite{wk1} under suitable settings for the classification
constraints. In order to find out which of the calculi listed there can be
restricted to~$\CPq{N-1}$, one tries to calculate 
the bimodule structures of restricted calculi. To this purpose,
expressions of the type $\d x\cdot y$ with $x,y\in\CPq{N-1}$ need to be
transformed into left-module expressions using the bimodule structure of
a $\Sq^{2N-1}$-calculus. The crucial question then is
how to recognise which expressions
on the right-hand side can be written in terms of $\CPq{N-1}$ only.
Moreover, the left-module relations which may hold in the module of
one-forms over $\CPq{N-1}$ have to be described.
Up to this point, we have no effective algorithm to solve these two problems.

The second approach consists in direct classification of covariant
first order differential calculi on $\CPq{N-1}$ under appropriate
algebraic constraints, using representation theory similarly as done
for Podle{\'s}' quantum spheres in~\cite{ap} or for the Vaksman-Soibelman
quantum spheres in~\cite{wk1}. However, choosing appropriate constraints
turns out more difficult than in the case of, e.\,g., the quantum spheres
$\Sq^{2N-1}$. One frequently used constraint setting requires the
differentials of the algebra generators to generate the bimodule of
one-forms as a free left module. Unlike for many other examples,
this standard setting is obviously inadequate here from a geometrical point
of view since the dimension of the differential calculus would then be
much higher than that of the algebra itself. Although we are going to consider
this setting, we shall look for more appropriate constraints. These should
at one hand be some kind of a natural choice while on the other hand they
should reduce the dimension of the differential calculus
such that it becomes close to that of the underlying quantum space.

Fortunately, a combination of the two approaches makes it much easier to
overcome the difficulties in both of them. The consideration of
co-representations which is the first step in following the classification
strategy, leads to a precise knowledge about the types of expressions that
may occur in the bimodule structure of covariant first order differential
calculi on $\CPq{N-1}$. Thus, it is much easier to decide whether
restricting a given $\Sq^{2N-1}$ calculus leads in fact to a
differential calculus on the sub-algebra.
On the other hand, if some differential calculi on $\CPq{N-1}$ can be
obtained by restricting calculi from the quantum spheres, the relations
which hold in these calculi will give evidence which type of
algebraic relations should be introduced into the classification constraint.

Nevertheless, the proof of the classification results still involves
rather complicated calculations which require the aid of computer algebra.

\subsection{First result: restriction of calculi from the quantum spheres}
\label{s-restr}
In~\cite{wk1} differential calculi on the quantum spheres~$\Sq^{2N-1}$
were classified. Two different classification constraints were applied.
(The classification is complete for $N\ge4$ but all calculi
described exist for $N\ge2$, too.)
We recall the main results: There are four families 
$\Gammaat$, $\Gammaap$, $\Gammap$, $\Gammapp$ of
covariant first order differential $*$\nob-calculi
with $\{\d z_1,\dots,\d z_N,\d z^*_1,\dots,\d z^*_N\}$ as a free left module
basis for the bimodule of one-forms. Each has two real parameters, with the
exception of certain parameter pairs for $\Gammaap$ with non-real $\alpha$
which we shall not include in our consideration. Further, there are three
families $\Gammatl$, $\Gammatp$, $\Gammatpp$ of covariant first order
differential $*$\nob-calculi, each with one real parameter, for which the
bimodule of one-forms is generated as a left module by $\{\d z_1,\dots,\d z_N,
\d z^*_1,\dots,\d z^*_N\}$ and for which all algebraic relations in the
left module of one-forms are generated by one relation
$\sumR{i}z_i\d z^*_i+\lambda\sumL{i}z^*_i\d z_i=0$.
(We restrict ourselves to $*$\nob-calculi here, leaving aside $\Gammato$ and
$\Gammatoo$.)
The equations taken verbatim from~\cite{wk1} which characterise the
bimodule structure of these calculi, are given in section~\ref{sRest},
equations~(\ref{eGatB})\nob--(\ref{eGtpp8B}).
We can now state our first result.
\begin{Thm}\label{tRest}
Each of the covariant first order differential $*$-calculi
$\Gammaat$, $\Gammaap$, $\Gammap$, $\Gammapp$, $\Gammatl$, $\Gammatp$,
$\Gammatpp$ on~$\Sq^{2N-1}$ can be restricted to a covariant first order
differential calculus on~$\CPq{N-1}$. All of the calculi $\Gammaat$,
$\Gammapp$ and $\Gammatpp$ yield the same restricted calculus
$\GCPt$ while all of the calculi $\Gammaap$, $\Gammap$, $\Gammatl$
and $\Gammatp$ lead to the same restricted calculus $\GCPtt$, independent
on the values of all parameters involved.

In $\GCPt$, all relations in the left module of one-forms are generated
by the set of relations ($1\le i,j\le N$)
\begin{align*}
\d\x ij&=q^2\sumL{s}\x is\d\x sj
+q^{-1}\sumR a\Rlm tubc \Rm sbia \Rc cvaj\x st\d\x uv
-\frac{\sii+}{\si+}\x ij\HCP\\
&\quad\qquad-\frac1{\si+}\delta_{ij}q^{2j}\HCP.
\end{align*}
The bimodule structure of $\GCPt$ is given by
\begin{align*}
\d\x ij\x kl
&=q^{-2}\sumR{a}\RCPm xyzwijkl\Rlm tubc \Rm sbza \Rc cvaw\x xy\x st\d\x uv\\
&\quad\qquad+q^3\RCPc stuvijkl\sumL{w}\x st\x uw\d\x wv
-\frac{q^{2N+2}}{\si+}\delta_{jk}\x ij\HCP\\
&\quad\qquad-\frac{q^{-1}}{\si+}\Rr abjk \Rm scia \Rc cvbl\x sv\HCP
-\frac{q^{-2}}{\si+}\delta_{ij}q^{2j}\x kl\HCP
-\frac{q^{-2}\siv+}{\si+}\x ij\x kl\HCP.
\end{align*}

In $\GCPtt$, all relations in the left modules of one-forms are generated
by the set of relations ($1\le i,j\le N$)
\begin{align*}
\d\x ij&=q^2\sumL{s}\x is\d\x sj
+q^{-1}\sumR a\Rlm tubc \Rm sbia \Rc cvaj\x st\d\x uv;& 
\HCP&=0.
\end{align*}
The bimodule structure of $\GCPtt$ is given by
\begin{align*}
\d\x ij\x kl
&=q^{-2}\sumR{a}\RCPm xyzwijkl\Rlm tubc \Rm sbza \Rc cvaw\x xy\x st\d\x uv\\
&\quad\qquad+q^3\RCPc stuvijkl\sumL{w}\x st\x uw\d\x wv.
\end{align*}
\end{Thm}

\noindent\emph{Remarks:}
1.\ Note that the first group of relations and the bimodule
structure of~$\GCPtt$ are obtained from the
relations and bimodule structure of~$\GCPt$
just by inserting the additional relation~$\HCP=0$. Thus,
the calculus $\GCPtt$ is obtained by factorising~$\GCPt$ by this relation.\\
2.\ It was proved in~\cite{sch} that the calculi $\Gammatl$ and
$\Gammatpp[q^{2N+2}]$ on $\Sq^{2N-1}$ are restrictions of covariant
(as for $\Gammatpp[q^{2N+2}]$, even bicovariant) first order differential
calculi on $\SUq(N)$. Since restriction of $\Gammatl$ to $\CPq{N-1}$ yields
$\GCPtt$, while $\Gammatpp[q^{2N+2}]$ can be restricted to $\GCPt$,
both $\GCPt$ and $\GCPtt$ are even restrictions of covariant differential
calculi on~$\SUq(N)$.

\subsection{Second result: direct classification}
For a classification of covariant first order differential calculi on
the quantum projective space $\CPq{N-1}$, we need a plausible constraint
setting which should essentially consist of a dimension restriction for the
bimodule of one-forms.
We shall consider three settings for the classification constraint.
First, we require that the bimodule of one-forms be generated
by $\d\x ij$, $1\le i,j\le N$, $(i,j)\ne(N,N)$ as free left module basis.
As stated above, this constraint does not make much sense from the
geometrical point of view since it means that the module of one-forms
needs to be of far higher dimension than the underlying algebra itself.
We consider this setting mostly for algebraic completeness since this
type of condition is the starting-point for all other types of dimension
condition taken into consideration.

The other two settings---which are supposed to be of geometrical relevance---%
are motivated by Theorem~\ref{tRest}. We suppose that the differential
calculi obtained by restriction of calculi from~$\Sq^{2N-1}$ have 
appropriate dimension. Therefore, we choose the relations found in
the calculi $\GCPt$ and $\GCPtt$ as ``templates'' for our second and
third constraint setting; the actual constraints are obtained by allowing
the coefficients in the relations to vary.

It turns out that allowing the coefficients to vary is in fact no essential
generalisation because the values taken by the coefficients in the relations
of the restricted calculi $\GCPt$ and $\GCPtt$ are the only possible ones;
finally, covariant differential calculi exist and are uniquely determined
under all three constraints.
The following theorem states our classification results for all these
settings.
\begin{Thm}\label{tClass}~
\begin{enumerate}
\item\label{it-free}
There is a covariant first order differential calculus
$\Gammad=\GCPd$ on $\CPq{N-1}$ for which
$\{\d\x ij~\vert~i,j=1,\dots,N;(i,j)\ne(N,N)\}$ is a free left module
basis of~$\Gamma$. If $N\ge6$, then $\GCPd$ is the only differential
calculus with this property. The bimodule structure of~$\GCPd$ is given by
\begin{align*}
\d\x ij\cdot\x kl
&=q^{-1}\RCPm stuvijkl\x st\d\x uv +q\RCPc stuvijkl\x st\d\x uv
+\RCPm stuvwxyz\RCPc wxyzijkl\x st\d\x uv\\
&\quad\qquad
-\sumL{w}\x ij\x kw\d\x wl -\sumL{w}\RCPc stuvijkl\x st\x uw\d\x wv\\
&\quad\qquad-q \sumR{a} \Rlm tubc \Rm sbka \Rc cval \x ij\x st\d\x uv\\
&\quad\qquad
-\sumR{a}\RCPm xyzwijkl\Rlm tubc \Rm sbza \Rc cvaw\x xy\x st\d\x uv
+(q^2+1)\x ij\x kl\HCP
\end{align*}
\item\label{it-red1}
For $N\ge6$, there is exactly one covariant first order differential calculus
$\Gammad$ on $\CPq{N-1}$ for which
$\{\d\x ij~\vert~i,j=1,\dots,N;(i,j)\ne(N,N)\}$
generates $\Gamma$ as a left module, and for which all relations in
the left module $\Gamma$ are algebraically generated by the set of relations
($1\le i,j\le N$)
\begin{align*}
\d\x ij
&=A\sumL{s}\x is\d\x sj+B\sumR a\Rlm tubc \Rm sbia \Rc cvaj\x st\d\x uv\\
&\quad\qquad+C\x ij\HCP+D\delta_{ij}q^{2j}\HCP
\end{align*}
for some fixed coefficients $A$, $B$, $C$, and $D$. This is the
calculus $\Gammad=\GCPtd$ from Theorem~\ref{tRest} with
\begin{align*}
A&=q^2;&B&=q^{-1};&C&=-\frac{\sii+}{\si+};&D&=-\frac1{\si+}.
\end{align*}
\item\label{it-red2}
For $N\ge6$, there is exactly one covariant first order differential calculus
$\Gammad$ on $\CPq{N-1}$ for which
$\{\d\x ij~\vert~i,j=1,\dots,N;(i,j)\ne(N,N)\}$
generates $\Gamma$ as a left module, and for which all relations in
the left module $\Gamma$ are algebraically generated by the set of relations
($1\le i,j\le N$)
\begin{equation*}
\d\x ij
=A\sumL{s}\x is\d\x sj+B\sumR a\Rlm tubc \Rm sbia \Rc cvaj\x st\d\x uv;\qquad
\HCP=0
\end{equation*}
for some fixed coefficients $A$, $B$. This is the
calculus $\Gammad=\GCPttd$ from Theorem~\ref{tRest} with
\begin{align*}
A&=q^2;&B&=q^{-1}.
\end{align*}

\end{enumerate}
\end{Thm}
\noindent\emph{Remarks:} 1.\
It is easily seen that, like $\GCPtt$ from $\GCPtt$, even $\GCPt$ (and thus,
$\GCPtt$) is obtained from $\GCP$ by factorisation because
the bimodule structure of $\GCP$ turns into that of $\GCPt$
if simply the 
set of left-module relations of $\GCPt$ is imposed.\\
2.\ All differential calculi discussed exist for $N\ge2$. For the
second and third case this is part of the statement of Theorem~\ref{tRest},
so it had to be stated explicitly here only for case~\ref{it-free}. However,
the uniqueness is guaranteed only for $N\ge6$; the reasons will become clear
from the proof (see sections~\ref{s-rt} and~\ref{sClass}).

\section{Representation theory}\label{s-rt}

We start by investigating corepresentations of the quantum group $\SUq(N)$
on $\CPq{N-1}$. By the coaction $\DeltaR$, a corepresentation
of $\SUq(N)$ on $\CPq{N-1}$ is given which decomposes into
summands corresponding to invariant vector spaces $V(k)$,
$k=0,1,\dots$ Here, $V(k)$ is the vector space of homogeneous
polynomials which is generated by
precisely those monomials of degree $k$ in the generators $\x ij$ which
are not reduced to lower degree by the algebra relations of $\CPq{N-1}$.
We denote by $\pi(k)$ the corepresentation of $\SUq(N)$ on $V(k)$.

Since $q$ is not a root of unity, the representation theory is essentially
identical to the classical case---see \cite{hay}---and the decomposition
of corepresentations into irreducible summands can be described
by means of Young frames, cf.\ \cite{br}. We shall use this notation in
the following, denoting the trivial corepresentation by~$\ynull$.

The bimodule structure of any covariant first order differential calculus
needs to be formed by intertwining morphisms $T\in\mathrm{Mor}((\pi(1)+\pi(0))
\otimes(\pi(1)+\pi(0)),\pi(k)\otimes(\pi(1)+\pi(0)))$. (Note that the
vector space generated by $\x ij$, $1\le i,j\le N$, is $V(0)\oplus V(1)$.)

We show the calculations with Young frames for $N=5$. For other $N\ge4$,
the calculations are analogous, while for $N=2$, $N=3$ additional
coincidences of Young frames and, thus, irreducible summands have to
be observed, spoiling the uniqueness argument in these cases. 
We shall have to sharpen the requirement even to $N\ge6$ because of
a linear independence argument used later.

Starting from $\pi(0)= \ynull$ and $\pi(1)=\ytt211141000$, one calculates
successively $\pi(k)$ and $\pi(k)\otimes(\pi(1)+\pi(0))$, $k=0,1,\dots$; note
that $\pi(k+1)$ is obtained from $\pi(k)$ by cancelling all those summands
from $\pi(k)\otimes(\pi(1)+\pi(0))$ which correspond to invariant subspaces
annihilated by the algebra relations of $\CPq{N-1}$.
\begin{align*}
\pi(0)&= \ynull,& \pi(1)&=\ytt211141000,\\
\pi(2)&=\ytt422244110,& \pi(3)&=~\yttv533344411\yttv100010000~,\quad\dots;
\end{align*}
\begin{equation*}
(\pi(1)+\pi(0))\otimes(\pi(1)+\pi(0))
=\ytt422244110 + \ytt311031100 + \ytt332244200 + \ytt221032000
    + 4\ytt211141000 + 2\ynull;
\end{equation*}
\begin{align*}
\pi(0)\otimes(\pi(1)+\pi(0))
&=\ytt211141000+\ynull,\\
\pi(1)\otimes(\pi(1)+\pi(0))
&=\ytt422244110 + \ytt311031100 + \ytt332244200 + \ytt221032000
    + 3\ytt211141000+\ynull,\\
\pi(2)\otimes(\pi(1)+\pi(0))&=~\yttv533344411\yttv100010000~ 
    + \ytt522143111 + \ytt543344421
    + \ytt432143210+ 3\ytt422244110\\
   &\quad\qquad + \ytt311031100
    + \ytt332244200
    + \ytt211141000,\\
\pi(3)\otimes(\pi(1)+\pi(0))&=~\yttv544444441\yttv300011100~
    + ~\yttv533244311\yttv200011000~
    + ~\yttv554444442\yttv200011000~ + ~\yttv543244321\yttv100010000~
    + 2~\yttv533344411\yttv100010000~\\
   &\quad\qquad + \ytt522143111
    + \ytt543344421 + \ytt422244110
\end{align*}%
By comparing these decompositions, it is found that
$\pi(k)\otimes(\pi(1)+\pi(0))$ has two summands in common with
$(\pi(1)+\pi(0))\otimes(\pi(1)+\pi(0))$ for $k=0$, five for $k=1$,
four for $k=2$, one for $k=3$, and none for higher $k$.
Taking into account the multiplicities
of all these summands---e.g.\ $\ytt211141000$ occurs with multiplicity
$4$ in $(\pi(1)+\pi(0))\otimes(\pi(1)+\pi(0))$, and $5$ in
$\sum\limits_k \pi(k)\otimes(\pi(1)+\pi(0))$, allowing for $20$
independent subspace mappings---it can be seen that up to $33$ morphisms
can occur. Because of the necessary condition
$\sum\limits_{i=1}^{N}\d\x ii=0$ some of the summands vanish automatically,
and we are left with a general ansatz for the bimodule structure of
a covariant differential calculus on $\CPq{N-1}$ containing $27$
morphisms, namely
\begin{align}
\kern2em&\kern-2em\d\x ij\x kl=
   a_1\x ij\d\x kl
 + a_2\RCPm stuvijkl\x st\d\x uv
 + a_3\RCPc stuvijkl\x st\d\x uv
\notag\\\phantom{\smallsumL a}&
 + a_4(\RCP{}{}{}{}{}{}{}{}\RCPc{}{}{}{}{}{}{}{})^{stuv}_{ijkl}\x st\d\x uv
 + a_5\delta_{jk}\smallsumL s\x is\d\x sl
\notag\\\phantom{\smallsumL a}&
 + a_6\delta_{jk}\smallsumR a\Rlm tubc \Rm sbia \Rc cval\x st\d\x uv
 + a_7\Rr abjk \Rm scia \Rc cvbl \smallsumL t\x st\d\x tv
\notag\\\phantom{\smallsumL a}&
 + a_8\smallsumR b\Rr stab \Rrm uvbc \Rl adij \Rlm dckl\x st\d\x uv
 + a_9\delta_{jk}\delta_{il}q^{-2l}\HCP
\notag\\\phantom{\smallsumL a}&
 + e_1\delta_{ij} q^{-2j}\smallsumL{s}\x ks\d\x sl
 + e_2\delta_{ij} q^{-2j}\smallsumR a\Rlm tubc \Rm sbka \Rc cval\x st\d\x uv
\notag\\\kern-2em\phantom{\smallsumL a}&
 + e_3\delta_{kl} q^{-2k}\smallsumL{s}\x is\d\x sj
 + e_4\delta_{kl} q^{-2k}\smallsumR a\Rlm tubc \Rm sbia \Rc cvaj\x st\d\x uv
\label{g-cpans}
\\\phantom{\smallsumL a}&
 + f_1\delta_{jk}\d\x il
 + f_2\Rr abjk \Rm scia \Rc cvbl\d\x sv
 + f_3\delta_{ij}q^{-2j}\d\x kl
\notag\\\phantom{\smallsumL a}&
 + f_4\delta_{kl}q^{-2k}\d\x ij
 + f_5\delta_{ij}\delta_{kl}q^{-2j-2k}\HCP
\notag\\\phantom{\smallsumL a}&
 + b_1\smallsumL{s}\x ij\x ks\d\x sl
 + b_2\RCPc stuvijkl\smallsumL{w}\x st\x uw\d\x wv
\notag\\\phantom{\smallsumL a}&
 + b_3\smallsumR a\Rlm tubc \Rm sbka \Rc cval\x ij\x st\d\x uv
 + b_4\RCP xyzwijkl\Rlm tubc \Rm sbza \Rc cvaw\x xy\x st\d\x uv
\notag\\\phantom{\smallsumL a}&
 + b_5\delta_{jk}\x ij\HCP
 + b_6\Rr abjk \Rm scia \Rc cvbl\x sv\HCP
\notag\\\phantom{\smallsumL a}&
 + g_1\delta_{ij}q^{-2j}\x kl\HCP
 + g_2\delta_{kl}q^{-2k}\x ij\HCP
 + c\x ij\x kl\HCP.
\notag
\end{align}
In case of the ``free'' classification constraint, i.e.\ if $\d\x ij$,
$(i,j)\ne(N,N)$ are supposed to be a free left module basis for $\Gamma$,
all morphisms are independent. If left-module relations among the $\d\x ij$
are admitted, some summands become superfluous, resulting in an ansatz
with less than $27$ coefficients.
\section{Proof of the theorem on restricted calculi}\label{sRest}
In~\ref{s-restr}, we mentioned the classification results on covariant
first order differential $*$\nob-calculi from~\cite{wk1}. We want to
state first the systems of equations describing the bimodule structure
of these differential calculi. Note that we use the abbreviations
$\H+:=\sumR i z_i\d z^*_i$ and $\H-:=\sumL i z^*_i\d z_i$ for the two basic
invariant one-forms. In the first four families, the $\d z_i$ and $\d z^*_i$,
$i=1,\dots,N$, form a free left-module basis for~$\Gamma$.
\begin{align}
&\begin{array}{@{}r@{~}l@{~}l@{}}
\hbox to1cm{\kern-8mm$\Gammaat:$\hfil
$\d z_k z_l$}&=q\alpha\Rm stkl z_s\d z_t&{}+(q^2\alpha-1)z_k\d z_l\\
  &&{}+q^2\alpha^2(1-\si+\tau)z_kz_l\H+\\
  &&{}+q^2(1-\alpha\si+\tau)z_kz_l\H-
  \\[1ex]
\d z^*_k z^*_l&=
  q^{-1}\alpha^{-1}\Rc stkl z^*_s\d z^*_t
  &{}+(q^{-2}\alpha^{-1}-1)z^*_k\d z^*_l\\
  &&{}+(1-\si+\tau)z^*_k z^*_l\H+\\
  &&{}+\alpha^{-2}(1-\alpha\si+\tau)z^*_kz^*_l\H-\\[1ex]
\d z_k z^*_l&=
  q^{-1}\alpha^{-1}\Rlm stkl z^*_s\d z_t&{}+(q^2\alpha-1)z_k\d z^*_l\\
  &&{}-q^2\alpha(1-\s+\tau)z_kz^*_l\H+  -\alpha\tau q^{2k}\delta_{kl}\H+\\
  &&{}-\alpha^{-1}(1-q^2\alpha\s+\tau)z_k z^*_l\H-  -\tau q^{2k}\delta_{kl}\H-
  \\[1ex]
\d z^*_k z_l&=
  q\alpha\Rr stkl z_s\d z^*_t&{}+(q^{-2}\alpha^{-1}-1)z^*_k\d z_l\\
  &&{}-q^2\alpha(1-\s+\tau)z^*_kz_l\H+  -q^{2N}\alpha\tau \delta_{kl}\H+\\
  &&{}-\alpha^{-1}(1-q^2\alpha\s+\tau)z^*_kz_l\H-  -q^{2N}\tau \delta_{kl}\H-
\end{array}
\label{eGatB}\\[2pt plus 2ex]
&\begin{array}{@{}r@{~}l@{~}l@{}}
\hbox to1cm{\kern-8mm$\Gammaap:$\hfil
$\d z_k z_l$}&=q\alpha\Rm stkl z_s\d z_t&{}+(q^2\alpha-1)z_k\d z_l\\
  &&{}+\omega z_kz_l\H++(\alpha^{-1}\omega
  -q^2(\alpha-1))z_kz_l\H-\\[1ex]
\d z^*_k z^*_l&=
  q^{-1}\alpha^{-1}\Rc stkl z^*_s\d z^*_t
  &{}+(q^{-2}\alpha^{-1}-1)z^*_k\d z^*_l\\
  &&{}+(q^2\alpha\omega^{-1}-(\alpha^{-1}-1))z^*_k z^*_l\H+\\
  &&{}+q^2\omega^{-1}z^*_kz^*_l\H-\\[1ex]
\d z_k z^*_l&=
  q^{-1}\alpha^{-1}\Rlm stkl z^*_s\d z_t&{}+(q^2\alpha-1)z_k\d z^*_l\\
  &&{}-q^2\alpha z_kz^*_l\H+-\alpha^{-1}z_k z^*_l\H-\\[1ex]
\d z^*_k z_l&=
  q\alpha\Rr stkl z_s\d z^*_t&{}+(q^{-2}\alpha^{-1}-1)z^*_k\d z_l\\
  &&{}-q^2\alpha z^*_kz_l\H+-\alpha^{-1}z^*_k z_l\H-
\end{array}
\label{eGapB}\\[2pt plus 2ex]
&\begin{array}{@{}r@{~}l@{~}l@{~}l@{}}
\hbox to1cm{\kern-8mm$\Gammap:$\hfil
$\d z_k z_l$}&=q^{-1}\Rm stkl z_s\d z_t
  &{}+\omega z_kz_l\H+&{}+(q^2\omega\psi-1)z_kz_l\H-\\[1ex]
\d z^*_k z^*_l&=
  q\Rc stkl z^*_s\d z^*_t
  &{}+(\psi-q^2)z^*_k z^*_l\H+ &{}+q^2\omega^{-1}z^*_kz^*_l\H-\\[1ex]
\d z_k z^*_l&= q\Rlm stkl z^*_s\d z_t &{}-z_kz^*_l\H+&{}-q^2z_k z^*_l\H-\\[1ex]
\d z^*_k z_l&= q^{-1}\Rr stkl z_s\d z^*_t&{}-z^*_kz_l\H+&{}-q^2z^*_k z_l\H-
\end{array}
\label{eGpB}\\[2pt plus 2ex]
&\begin{array}{@{}r@{~}l@{~}l@{}}
\hbox to1cm{\kern-8mm$\Gammapp:$\hfil
$\d z_k z_l$}&=q^{-1}\Rm stkl z_s\d z_t
  &{}-q^{-2}\dfrac{\varrho}{\tau}(\si+\varrho-1)z_kz_l\H+
  -\dfrac{\varrho}{\tau}(\si+\tau-q^2)z_kz_l\H-
\kern-7mm
\\[1ex]
\d z^*_k z^*_l&=
  q\Rc stkl z^*_s\d z^*_t
  &{}-\dfrac{\tau}{\varrho}(\si+\varrho-1)z^*_k z^*_l\H+
  -q^2\dfrac{\tau}{\varrho}(\si+\tau-q^2)z^*_kz^*_l\H-
\kern-7mm
\\[1ex]
\d z_k z^*_l&= q\Rlm stkl z^*_s\d z_t
  &{}-q^{-2}\varrho q^{2k}\delta_{kl}\H+-\tau q^{2k}\delta_{kl}\H-\\
  &&{}+(\s+\varrho-1)z_kz^*_l\H++q^2(\s+\tau-1)z_kz^*_l\H-\\[1ex]
\d z^*_k z_l&= q^{-1}\Rr stkl z_s\d z^*_t
  &{}-q^{2N-2}\varrho\delta_{kl}\H+-q^{2N}\tau\delta_{kl}\H-\\
  &&{}+(\s+\varrho-1)z^*_kz_l\H++q^2(\s+\tau-1)z^*_kz_l\H-
\end{array}
\label{eGppB}
\end{align}
In the following three families of calculi, $\d z_i$, $\d z^*_i$ 
still generate $\Gamma$ as a left module but no longer as a free one.
Instead, all left-module relations are algebraically generated
by $\H++\lambda\H-=0$ where $\lambda$ is the (real) parameter of the families
of calculi.
\begin{align}
&\begin{array}{@{}r@{~}l@{~}l@{~}l@{}}
\hbox to1cm{\kern-8mm$\Gammatl:$\hfil
$\d z_k z_l$}&=q\lambda^{-1}\Rm stkl z_s\d z_t&{}+(q^2\lambda^{-1}-1)z_k\d z_l
  &{}+q^2\lambda^{-1}(\lambda^{-1}-1)z_kz_l\H+\\[1ex]
\d z^*_k z^*_l &=q^{-1}\lambda\Rc stkl z^*_s\d z^*_t
  &{}+(q^{-2}\lambda-1)z^*_k\d z^*_l
  &{}-(\lambda-1)z^*_kz^*_l\H+\\[1ex]
\d z_k z^*_l &=q^{-1}\lambda\Rlm stkl z^*_s\d z_t
  &{}+(q^2\lambda^{-1}-1)z_k\d z^*_l
  &{}-(q^2\lambda^{-1}-1)z_kz^*_l\H+\\[1ex]
\d z^*_k z_l &=q\lambda^{-1}\Rr stkl z_s\d z^*_t
  &{}+(q^{-2}\lambda-1)z^*_k\d z_l
  &{}-(q^2\lambda^{-1}-1)z^*_kz_l\H+
\end{array}
\label{eGtlB}\\[2pt plus 2ex]
&\begin{array}{@{}r@{~}l@{~}l@{}}
\hbox to1cm{\kern-8mm$\Gammatp:$\hfil
$\d z_k z_l$}&=q^{-1}\Rm stkl z_s\d z_t
  &{}-\lambda^{-1}(q^4\lambda^{-1}-1)z_kz_l\H+\\[1ex]
\d z^*_k z^*_l &=q\Rc stkl z^*_s\d z^*_t
  &{}-q^{-2}\lambda(q^4\lambda^{-1}-1)z^*_kz^*_l\H+\\[1ex]
\d z_k z^*_l &=q\Rlm stkl z^*_s\d z_t
  &{}+(q^2\lambda^{-1}-1)z_kz^*_l\H+\\[1ex]
\d z^*_k z_l &=q^{-1}\Rr stkl z_s\d z^*_t
  &{}+(q^2\lambda^{-1}-1)z^*_kz_l\H+
\end{array}
\label{eGtpB}\\[2pt plus 2ex]
&\begin{array}{@{}r@{~}l@{~}l@{}}
\hbox to1cm{\kern-8mm$\Gammatpp,
\lambda\not\in\{0,\infty\}:$\kern-20mm\hfil}&\\[1ex]
\d z_k z_l&=q^{-1}\Rm stkl z_s\d z_t\\[1ex]
\d z^*_k z^*_l &=q\Rc stkl z^*_s\d z^*_t\\[1ex]
\d z_k z^*_l &=q\Rlm stkl z^*_s\d z_t
  &{}+q^{-2}\si+^{-1}(q^4\lambda^{-1}-1)q^{2k}\delta_{kl}\H+\\
  &&{}-\si+^{-1}(q^{2N+2}\lambda^{-1}-1)z_kz^*_l\H+\\[1ex]
\d z^*_k z_l &=q^{-1}\Rr stkl z_s\d z^*_t
  &{}+q^{2N-2}\si+^{-1}(q^4\lambda^{-1}-1)\delta_{kl}\H+\\
  &&{}-\si+^{-1}(q^{2N+2}\lambda^{-1}-1)z^*_kz_l\H+
\end{array}
\label{eGtppB}\\[2pt plus 2ex]
&\begin{array}{@{}r@{~}l@{}}
\hbox to1cm{\kern-8mm$\Gammatpp[0]:$\hfil
$\d z_k z_l$}&=q^{-1}\Rm stkl z_s\d z_t\\[1ex]
\d z^*_k z^*_l &=q\Rc stkl z^*_s\d z^*_t\\[1ex]
\d z_k z^*_l &=q\Rlm stkl z^*_s\d z_t
  -q^{-2N+2}\si+^{-1}q^{2k}\delta_{kl}\H-
  +q^2\si+^{-1}z_kz^*_l\H-\\[1ex]
\d z^*_k z_l &=q^{-1}\Rr stkl z_s\d z^*_t
  -q^2\si+^{-1}\delta_{kl}\H-
  +q^2\si+^{-1}z^*_kz_l\H-
\end{array}
\label{eGtpp0B}\\[2pt plus 2ex]
&\begin{array}{@{}r@{~}l@{}}
\hbox to1cm{\kern-8mm$\Gammatpp[\infty]:$\hfil
$\d z_k z_l$}&=q^{-1}\Rm stkl z_s\d z_t\\[1ex]
\d z^*_k z^*_l &=q\Rc stkl z^*_s\d z^*_t\\[1ex]
\d z_k z^*_l &=q\Rlm stkl z^*_s\d z_t
  -q^{-2}\si+^{-1}q^{2k}\delta_{kl}\H+
  +\si+^{-1}z_kz^*_l\H+\\[1ex]
\d z^*_k z_l &=q^{-1}\Rr stkl z_s\d z^*_t
  -q^{2N-2}\si+^{-1}\delta_{kl}\H+
  +\si+^{-1}z^*_kz_l\H+
\end{array}
\label{eGtpp8B}
\end{align}%
To compute restrictions of these calculi to the sub-algebra $\CPq{N-1}$, 
we exploit the results of the preceding section. 
We rewrite $\d \x ij\cdot\x kl$ and all the $27$ summands occurring on the
right-hand side of~\eqref{g-cpans} in terms of the generators of $\Sq^{2N-1}$
via the embedding $\x ij=z_iz^*_j$ and Leibniz rule. After transforming
all these terms to left-module expressions by applying the bimodule
structure of a given differential calculus over~$\Sq^{2N-1}$, 
the ansatz equation~\eqref{g-cpans} is written in left-module expressions
from the $\Sq^{2N-1}$ calculus. By comparing coefficients (carefully observing
given relations), it is easily determined whether equation~\eqref{g-cpans}
can be satisfied for appropriate values of the coefficients $a_1$, $a_2$,
\dots, $c$.---%
Moreover, since relations in a calculus must cancel out invariant
subspaces for certain co-representations, one can find out which relations
occur in a restricted calculus also by rewriting expressions for morphisms
in terms of $\Sq^{2N-1}$ and checking whether some of them become zero.

We demonstrate the procedure for $\Gammaat$ with the bimodule
structure~\eqref{eGatB}. For abbreviation, let $\H\alpha:=\alpha\H++\H-$.
From the equalities
\begin{gather*}
\begin{aligned}
\d\x ij&=q^{-1}\alpha^{-1}\Rlm stij z^*_s\d z_t
+q^2\alpha z_i\d z^*_j-q^2\alpha(1-\s+\tau)z_iz^*_j\H+
\\&\qquad{}-\alpha\tau\delta_{ij}q^{2i}\H+
-\alpha^{-1}(1-q^2\alpha\s+\tau)z_iz^*_j\H-
-\tau\delta_{ij}q^{2j}\H-\end{aligned}\\
\sumL k\x ik\d\x kj=\alpha z_i\d z^*_j
-\alpha(1-\si+\tau)z_iz^*_j\H+ +\si+\tau z_iz^*_j\H-\\
\begin{aligned}{}&\sumR a\Rlm tubc \Rm sbia \Rc cvaj \x st\d\x uv\\
&\qquad\qquad=\alpha^{-1}\Rlm stij z^*_s\d z_t +q\alpha\si+\tau z_iz^*_j\H+
-q\alpha^{-1}(1-\alpha\si+\tau)z_iz^*_j\H-\end{aligned}\\
\HCP=\si+\tau\H\alpha
\end{gather*}
it can be seen that the relation
\begin{equation}
\label{g-gcptrel}
\begin{aligned}
\d\x ij&=q^2\sumL k\x ik\d\x kj
+q^{-1}\sumR a\Rlm tubc \Rm sbia \Rc cvaj \x st\d\x uv\\
&\quad\qquad-\frac{\sii+}{\si+}\x ij\HCP -\frac1{\si+}\delta_{ij}q^{2j}\HCP
\end{aligned}
\end{equation}
holds. This implies immediately that the morphisms corresponding with
the coefficients $a_1$, $a_2$, $a_3$, $a_4$, $f_1$, $f_2$, $f_3$, $f_4$
in equation~\eqref{g-cpans} are linearly dependent on the other ones and
can be omitted. We therefore continue by rewriting expressions from the
right-hand side of~\eqref{g-cpans}---i.e.\ the summands that can occur
in the bimodule structure of a covariant differential calculus
over~$\CPq{N-1}$---in terms of the generators of~$\Sq^{2N-1}$.
\begin{align*}
{}&\begin{aligned}
{}&\d\x ij\cdot\x kl=
q^{-1}\alpha^{-1}\Rlm uvbl\Rm tbak\Rlm saij z^*_sz_tz^*_u\d z_v\\&\qquad
+q^2\alpha\Rc uval \Rr tajk z_iz_tz^*_u\d z^*_v
+(q^4-q^{-2}\alpha^{-1})z_iz^*_jz_kz^*_l\H-\\&\qquad
-q^{-1}\tau\delta_{ab} q^{2a}\Rm saic \Rc btdl \Rr cdjk z_sz^*_t\H\alpha
-q^{2N+2}\tau\delta_{jk}z_iz^*_l\H\alpha\\&\qquad
-q^{-2}\tau\delta_{ij}q^{2j}z_kz^*_l\H\alpha
+(-q^4+\tau+q^2\tau+q^4\s+\tau)z_iz^*_jz_kz^*_l\H\alpha
\end{aligned}\\
{}&\begin{aligned}
{}&\RCP xyzwijkl\Rlm tubc \Rm sbza \Rc cvaw\x xy\x st\d\x uv
=q\alpha^{-1}\Rlm uvbl\Rm tbak\Rlm saij z^*_sz_tz^*_u\d z_v\\&\qquad
+\alpha\si+\tau z_iz^*_jz_kz^*_l\H+
-\alpha^{-1}(1-\alpha\si+\tau)z_iz^*_jz_kz^*_l\H-
\end{aligned}\\
{}&\begin{aligned}
{}&\RCPc stuvijkl\smallsumL{w}\x st\x uw\d\x wv
=q^{-1}\alpha\Rc uval \Rr tajk z_iz_tz^*_u\d z^*_v\\&\qquad
-q\alpha(1-\si+\tau)z_iz^*_jz_kz^*_l\H+
+q\si+\tau z_iz^*_jz_kz^*_l\H-
\end{aligned}
\end{align*}
which, together with the preceding equations, imply
\begin{equation}
\begin{aligned}
{}&\d\x ij\cdot\x kl
-q^{-2}\RCP xyzwijkl\Rlm tubc \Rm sbza \Rc cvaw\x xy\x st\d\x uv\\
&\qquad\qquad
-q^3\RCPc stuvijkl\smallsumL{w}\x st\x uw\d\x wv\\
{}&\quad=-{\si+}^{-1}\bigl( q^{-2} \siv+\x ij\x kl
+q^{-1}\delta_{ab} q^{2a}\Rm saic \Rc btdl \Rr cdjk \x st\\&\qquad\qquad\qquad
+q^{2N+2}\delta_{jk}\x il +q^{-2}\delta_{ij}q^{2j}\x kl \bigr)\HCP
\end{aligned}
\end{equation}
and therefore the bimodule structure of $\GCPt$ stated in Theorem~\ref{tRest}.
By rewriting all the morphisms from the right-hand side 
of~\eqref{g-cpans}, it is easily seen that they display no more linear
dependencies than those implied by the relation~\eqref{g-gcptrel}.
From these considerations, it is clear that $\GCPt$ is the restriction
of~$\Gammaat$ to the subalgebra~$\CPq{N-1}$.

Doing completely analogous calculations one proves the other restriction
statements of Theorem~\ref{tRest}.
\qed

\section{Proof of the classification theorem}\label{sClass}
Our proof of Theorem~\ref{tClass}
is organised in four parts. The first part deals with the 
possible coefficients in the classification constraints of the second and
third case while the other three ones are devoted each to one of the constraint
settings, showing that in each case there exists one uniquely determined
calculus.
\subsection{Determination of possible coefficients of the constraint relations}
First we show that if all left-module relations in a covariant
first order differential calculus are to be generated by one family of
relations of type
\begin{align*}
\d\x ij
&=A\sumL{s}\x is\d\x sj+B\sumR a\Rlm tubc \Rm sbia \Rc cvaj\x st\d\x uv
+C\x ij\HCP+D\delta_{ij}q^{2j}\HCP
\end{align*}
then the coefficients $A$, $B$, $C$, and $D$ have to take the values which
hold for $\GCPt$.
In fact, substituting $\d\x jk$ in $\sumL j\x ij\d\x jk$ by the above
relation yields
\begin{align*}
\sumL j\x ij\d\x jk&=A\sumL j\sumL s\x ij\x js\d\x sk\\*&\quad\qquad
+B\sumL j\sumR a\Rlm tubc \Rm sbja \Rc cvak\x ij\x st\d\x uv\\*&\quad\qquad
+C\sumL j\x ij\x jk\HCP + D\sumL j\delta_{jk}q^{2k}\HCP\\
&=q^{-2}A\sumL j\x ij\d\x jk + (q^{-1}B+q^{-2}C+D)\x ik\HCP
\end{align*}
which implies by coefficient comparison
\begin{align*}
A&=q^2;&q^{-1}B+q^{-2}C+D&=0.
\end{align*}
An analogous substitution for $\d\x uv$ in
$\sumR a\Rlm tubc \Rm sbia \Rc cvaj\x st\d\x uv$ leads to
\begin{equation*}
\begin{split}
\sumR a\Rlm tubc \Rm sbia \Rc cvaj\x st\d\x uv
=qB\sumR a\Rlm tubc \Rm sbia \Rc cvaj\x st\d\x uv\quad\\
+\bigl(q^{-1}A-(q^2-1)B+q^{-1}C+qD\bigr)\x ij\HCP
\end{split}
\end{equation*}
and therefore
\begin{align*}
B&=q^{-1};&q^{-1}A-(q^2-1)B+q^{-1}C+qD&=0.
\end{align*}
By inserting the values of $A$ and $B$, the last equation simplifies to
\begin{align*}
q^{-2}C+D&=q^{-2}.
\end{align*}
Finally, inserting the assumed relation in $\sumR i\d\x ii=0$ gives 
\begin{align*}
0&=A\sumR i\sumL j\x ij\d\x ji
+ B\sumR i\sumR a\Rlm tubc \Rm sbia \Rc cvai\x st\d\x uv\\&\quad\qquad
+ C\sumR i\x ii\HCP + D\sumR i q^{2i}\HCP\\
&=(A+qB+C+q^2\s+D)\HCP\\
\intertext{and thus}
0&=A+qB+C+q^2\s+D=q^2+q^2\si+D
\end{align*}
which implies 
\begin{align*}
D&=-\frac1{\si+};&C&=-\frac{\sii+}{\si+}.
\end{align*}

If all left-module relations in a covariant first order differential calculus
are assumed to be generated by the relations
\begin{align*}
\d\x ij
&=A\sumL{s}\x is\d\x sj+B\sumR a\Rlm tubc \Rm sbia \Rc cvaj\x st\d\x uv\\
\HCP&=0,
\end{align*}
similar (just easier) calculations as above lead directly to
\begin{align*}
A&=q^2;&B&=q^{-1}.
\end{align*}

Now we can prove the classification results for the different constraint
settings under consideration. To accomplish this, we use the respective
bimodule structure ansatzes and relations to evaluate necessary conditions
which result from the definition of differential calculus and the algebra
structure of~$\CPq{N-1}$.
Since it appears hopeless to do by hand the extensive calculations involved
(perhaps except for the case of the most reduced classification constraint),
a special-purpose computer algebra program written by the author was
employed to carry out the substitutions and term-reductions with
R-matrices. It should be emphasised that only substitutions of
R\nob-matrix expressions via given relations were done automatically,
a detailed discussion of the reduction strategy thus being not necessary.
The linear independence of the summands in the resulting expressions was
checked manually.
Because coefficients have to be compared for expressions containing
R\nob-matrices and $\x ij$, $\d\x ij$ with up to $6$ free indices, linear
independence is clear in some cases only for $N\ge6$; that's why this
assumption is made in the uniqueness statements of the theorem.

\subsection{Case~\ref{it-red2}}\label{s-red2}
By the assumed relations, the general ansatz~\eqref{g-cpans} is reduced
to 12 summands:
\begin{equation}
\kern-2em\begin{aligned}
\kern2em&\kern-2em\d\x ij\x kl=
   a_5\delta_{jk}\smallsumL s\x is\d\x sl
 + a_6\delta_{jk}\smallsumR a\Rlm tubc \Rm sbia \Rc cval\x st\d\x uv
\\\phantom{\smallsumL a}&
 + a_7\Rr abjk \Rm scia \Rc cvbl \smallsumL t\x st\d\x tv
 + a_8\smallsumR b\Rr stab \Rrm uvbc \Rl adij \Rlm dckl\x st\d\x uv
\\\phantom{\smallsumL a}&
 + e_1\delta_{ij} q^{-2j}\smallsumL{s}\x ks\d\x sl
 + e_2\delta_{ij} q^{-2j}\smallsumR a\Rlm tubc \Rm sbka \Rc cval\x st\d\x uv
\\\phantom{\smallsumL a}&
 + e_3\delta_{kl} q^{-2k}\smallsumL{s}\x is\d\x sj
 + e_4\delta_{kl} q^{-2k}\smallsumR a\Rlm tubc \Rm sbia \Rc cvaj\x st\d\x uv
\\\phantom{\smallsumL a}&
 + b_1\smallsumL{s}\x ij\x ks\d\x sl
 + b_2\RCPc stuvijkl\smallsumL{w}\x st\x uw\d\x wv
\\\phantom{\smallsumL a}&
 + b_3\smallsumR a\Rlm tubc \Rm sbka \Rc cval\x ij\x st\d\x uv
 + b_4\RCP xyzwijkl\Rlm tubc \Rm sbza \Rc cvaw\x xy\x st\d\x uv.
\end{aligned}\kern-2em
\label{g-cpansr2}
\end{equation}
We use the following necessary conditions for first order
differential calculi on~$\CPq{N-1}$ ($1\le i,j,k,l,m\le N$):
\begin{align}
\left(\sumR i\d\x ii\right)\x jk&=0 \label{kl1}\\
\d\x ij\sumR k\x kk-\d\x ij&=0 \label{kl2}\\
\sumL j\x ij\d\x jk-q^{-2}\d\x ik+\sumL j\d\x ij\cdot\x jk&=0 \label{kl3}\\
\d\x ij\cdot\x kl+\x ij\d\x kl
-q\RCPcm stuvijkl(\d\x st\cdot\x uv+\x st\d\x uv)
&=0 \label{kl4}\\
\d\x ij\cdot\x kl+\x ij\d\x kl
-q^{-1}\RCP stuvijkl(\d\x st\cdot\x uv+\x st\d\x uv)
&=0 \label{kl5}\\
\sumL m\d\x ij\cdot\x km\x ml-q^{-2}\d\x ij\cdot\x kl&=0 \label{kl7}
\end{align}
(the first two of which, as the last one, follow immediately from
$\sumR i\x ii=1$ while the
remaining ones are the result of deriving algebra relations via the
Leibniz rule).
Because of \eqref{g-cpansr2}, we obtain from condition~\eqref{kl1}
\begin{align*}
0&=(a_5+q^{2N+1}a_7+q^2\s+e_1+b_1)\sumL m\x jm\d\x mk\\*&\quad\qquad
+(a_6+q^{2N+1}a_8+q^2\s+e_2+b_3)\sumR a \Rlm tubc\Rm sbja\Rc cvak\x st\d\x uv.
\end{align*}
The two R\nob-matrix expressions on the right-hand side are linearly
independent such that we can compare coefficients to obtain two equations
that must hold for the coefficients of~\eqref{g-cpansr2}. The other conditions
are evaluated similarly.

By this method, the first five conditions provide us with a (redundant) set
of equations for the coefficients of~\eqref{g-cpansr2},
\begin{gather*}
\begin{aligned}
b_1+a_5+q^2\s+e_1+q^{2N+1}a_7&=0;\quad&
qb_3+qa_6+q^3\s+e_2+q^{2N+2}a_8&=0;\\
q^{-1}b_2+qa_7+a_5+q^2\s+e_3&=q^2;\quad&
q^2b_4+q^2a_8+qa_6+q^3\s+e_4&=1;\\
q^{-2}b_1+e_3+e_1+q^{-2N}\s+a_5&=0;&
b_4+qe_4+qe_2+q^{-2N+1}\s+a_6&=q^{-2};
\end{aligned}\\
\begin{aligned}
q^{-1}b_2&=b_1+q^2;\quad&e_1&=q^{-1}a_7;\quad&q^2e_3&=q^{-2N}a_5;\quad&
qe_2&=a_8;\\q^3e_4&=q^{-2N+1}a_6;\quad&b_4&=q^{-1}b_3+q^{-2};\quad&
q^{-2}e_3&=q^{-1}a_7;\quad&q^{-1}e_4&=a_8;
\end{aligned}\\
\begin{aligned}
e_1+(q-q^{-1})a_7-q^{-2N-2}a_5&=0;\quad&
qe_2+(q^2-1)a_8-q^{-2N-1}a_6&=0;
\end{aligned}\\
\begin{aligned}
a_5-q^{2N}e_1-q^{2N-1}(q^2-1)a_7-(1-q^{-2})a_5&=0;\\
qa_6-q^{2N-1}e_2-q^{2N}(q^2-1)a_8-(q-q^{-1})a_6&=0,
\end{aligned}
\end{gather*}
which is fulfilled if the coefficients depend on two complex parameters
$\alpha$, $\beta$ via
\begin{align*}
a_5&=\alpha;&a_6&=\beta;&a_7&=q^{-2N-3}\alpha;&a_8&=q^{-2N-3}\beta;\\
e_1&=q^{-2N-4}\alpha;&e_2&=q^{-2N-4}\beta;&e_3&=q^{-2N-2}\alpha;&
e_4&=q^{-2N-2}\beta;\\
b_1&=-(1+q^{-2}+q^{-2N-2}\s+)\alpha;\kern-8em&&&
b_2&=q^3-(q+q^{-1}+q^{-2N-1}\s+)\alpha;\kern-8em&&\\
b_3&=-(1+q^{-2}+q^{-2N-2}\s+)\beta;\kern-8em&&&
b_4&=q^{-2}-(q^{-1}+q^{-3}+q^{-2N-3}\s+)\beta.\kern-8em
\end{align*}
By evaluating condition~\eqref{kl7} using this two-parameter form of all
coefficients and doing coefficient comparison, we obtain finally
$\alpha=\beta=0$ which makes clear that $\GCPtt$ is the only covariant
first order differential calculus under this constraint setting.

\subsection{Case~\ref{it-red1}}\label{s-red1}
By the relations of this case, the bimodule structure ansatz~\eqref{g-cpans}
is reduced to 19 members, namely
\begin{equation}
\kern-2em\begin{aligned}
\kern2em&\kern-2em\d\x ij\x kl=
   a_5\delta_{jk}\smallsumL s\x is\d\x sl
 + a_6\delta_{jk}\smallsumR a\Rlm tubc \Rm sbia \Rc cval\x st\d\x uv
\\\phantom{\smallsumL a}&
 + a_7\Rr abjk \Rm scia \Rc cvbl \smallsumL t\x st\d\x tv
 + a_8\smallsumR b\Rr stab \Rrm uvbc \Rl adij \Rlm dckl\x st\d\x uv
\\\phantom{\smallsumL a}&
 + a_9\delta_{jk}\delta_{il}q^{-2l}\HCP
 + f_5\delta_{ij}\delta_{kl}q^{-2j-2k}\HCP
\\\phantom{\smallsumL a}&
 + e_1\delta_{ij} q^{-2j}\smallsumL{s}\x ks\d\x sl
 + e_2\delta_{ij} q^{-2j}\smallsumR a\Rlm tubc \Rm sbka \Rc cval\x st\d\x uv
\\\phantom{\smallsumL a}&
 + e_3\delta_{kl} q^{-2k}\smallsumL{s}\x is\d\x sj
 + e_4\delta_{kl} q^{-2k}\smallsumR a\Rlm tubc \Rm sbia \Rc cvaj\x st\d\x uv
\\\phantom{\smallsumL a}&
 + b_1\smallsumL{s}\x ij\x ks\d\x sl
 + b_2\RCPc stuvijkl\smallsumL{w}\x st\x uw\d\x wv
\\\phantom{\smallsumL a}&
 + b_3\smallsumR a\Rlm tubc \Rm sbka \Rc cval\x ij\x st\d\x uv
 + b_4\RCP xyzwijkl\Rlm tubc \Rm sbza \Rc cvaw\x xy\x st\d\x uv
\\\phantom{\smallsumL a}&
 + b_5\delta_{jk}\x ij\HCP
 + b_6\Rr abjk \Rm scia \Rc cvbl\x sv\HCP
\\\phantom{\smallsumL a}&
 + g_1\delta_{ij}q^{-2j}\x kl\HCP
 + g_2\delta_{kl}q^{-2k}\x ij\HCP
 + c\x ij\x kl\HCP.
\end{aligned}\kern-2em
\label{g-cpansr1}
\end{equation}
Now we evaluate again the conditions~\eqref{kl1}--\eqref{kl5}. By carrying
out the coefficient comparison completely (for all morphisms appearing) for
\eqref{kl4} and~\eqref{kl5} but only in part for the first three conditions
(since some of the expressions are very lengthy),
we obtain the equations
\begin{gather*}
\begin{aligned}
b_1+a_5+q^2\s+e_1+q^{2N+1}a_7&=0;\quad&
qb_3+qa_6+q^3\s+e_2+q^{2N+2}a_8&=0;\\
q^{-1}b_2+qa_7+a_5+q^2\s+e_3&=q^2;\quad&
q^2b_4+q^2a_8+qa_6+q^3\s+e_4&=1;\\
q^{-2}b_1+e_3+e_1-\s+a_5&=0;\quad&
b_4+qe_4+qe_2-q\s+a_6&=q^{-2};
\end{aligned}\\
\begin{aligned}
q^{-1}b_2&=b_1+q^2;\quad&
e_1&=q^{-1}a_7;\quad&
q^2e_3&=q^{-2N}a_5;\\
qe_2&=a_8;\quad&
q^3e_4&=q^{-2N+1}a_6;\quad&
g_1&=q^{-1}b_6;\\
b_4&=q^{-1}b_3+q^{-2};\quad&
q^{-2}e_3&=q^{-1}a_7;\quad&
q^{-1}e_4&=a_8;
\end{aligned}\\
\begin{aligned}
q^2g_2&=q^{-2N}b_5+q^2{\si+}^{-1};\quad&
e_1&=q^{-2N-2}a_5-(q-q^{-1})a_7;\\
qe_2&=q^{-2N-1}a_6-(q^2-1)a_8;\quad&
g_1&=q^{-2N-2}b_5-(q-q^{-1})b_6;\\
q^{-2}g_2&=q^{-1}b_6+q^{-2}{\si+}^{-1};\quad&
a_9&=-q^{2N+2}f_5
\end{aligned}
\end{gather*}
which reduces the number of independent coefficients to $7$, namely
$a_7$, $a_8$, $a_9$, $b_1$, $b_3$, $b_6$, and $c$, via the relations
\begin{align*}
a_5&=q^{2N+3}a_7;\quad&a_6&=q^{2N+3}a_8;\quad&
b_2&=qb_1+q^3;\quad&b_4&=q^{-1}b_3+q^{-2};\quad\\
e_1&=q^{-1}a_7;\quad&e_2&=q^{-1}a_8;\quad&
e_3&=qa_7;\quad&e_4&=qa_8;\\
b_5&=q^{2N+3}b_6;\quad&
f_5&=q^{-2N-2}a_9;\quad&
g_1&=q^{-1}b_6;\quad&g_2&=qb_6+{\si+}^{-1}.
\end{align*}
Further information is obtained from the condition
\begin{equation}
q^2\sumL m\x im\d\x mj\cdot\x kl
+q^{-1}\sumR a\Rlm tubc\Rm sbia\Rc cvaj\x st\d\x uv\cdot\x kl
-\d\x ij\cdot\x kl=0
\label{kl8}
\end{equation}
which is, of course, specific to this particular constraint setting
because it is derived from the imposed relation.
This equation leads to $a_7=a_8=b_1=b_3=a_9=0$ leaving just two parameters
$b_6$ and $c$.
With these simplification, we do the remaining coefficient comparisons
for condition~\eqref{kl1} and obtain
\begin{align*}
g_2&=0;\quad&c&=-q^{-2}-q^{-4}-(q\s++q^{2N+1}+q^{2N+3})b_6,
\end{align*}
therefore finally as the unique solution the coefficients of~$\GCPt$.

\subsection{Case~\ref{it-free} (free left module)}\label{s-free}
This is the most difficult case to handle. Here, the full
ansatz~\eqref{g-cpans} with $27$ unknown coefficients applies.
We start by evaluating the conditions \eqref{kl4} and \eqref{kl5}.
Coefficient comparison leads to the following equations for the coefficients
of~\eqref{g-cpans}:
\begin{gather*}
\begin{aligned}
a_3&=qa_1+q;\quad&
a_4&=qa_2;\quad&
f_3&=q^{-1}f_2;\quad&
f_4&=q^{-2N-2}f_1;\\
b_2&=qb_1;\quad&
e_1&=q^{-1}a_7;\quad&
e_3&=q^{-2N-2}a_5;\quad&
e_2&=q^{-1}a_8;\\
e_4&=q^{-2N-2}a_6;\quad&
g_1&=q^{-1}b_6;\quad&
g_2&=q^{-2N-2}b_5;\quad&
f_5&=q^{-2N-2}a_9;\\
a_2&=q^{-1}a_1+q^{-1};\quad&
a_4&=q^{-1}a_3;\quad&
f_4&=qf_2;\quad&
b_4&=q^{-1}b_3;\\
e_3&=qa_7;\quad&
e_4&=qa_8;\quad&
g_2&=qb_6;\quad&&
\end{aligned}\\
\begin{aligned}
f_3&=q^{-2N-2}f_1-(q-q^{-1})f_2;\quad&
e_1&=q^{-2N-2}a_5-(q-q^{-1})a_7;\\
e_2&=q^{-2N-2}a_6-(q-q^{-1})a_8;\quad&
g_1&=q^{-2N-2}b_5-(q-q^{-1})b_6.
\end{aligned}
\end{gather*}
Using these equations, we can rewrite the coefficients of the ansatz as
dependent on only $9$ parameters $\alpha$, $\beta$, $\gamma$, $\delta_1$,
$\delta_2$, $\delta_3$, $\varepsilon$, $\zeta$, and $c$ via
\begin{equation}\label{p9}
\begin{aligned}
a_1&=\alpha-1;&
a_2&=q^{-1}\alpha;&
a_3&=q\alpha;&
a_4&=\alpha;\\
a_5&=\beta;&
a_6&=\gamma;&
a_7&=q^{-2N-3}\beta;&
a_8&=q^{-2N-3}\gamma;\\
e_1&=q^{-2N-4}\beta;&
e_2&=q^{-2N-4}\gamma;&
e_3&=q^{-2N-2}\beta;&
e_4&=q^{-2N-2}\gamma;\\
b_1&=\delta_2;&
b_2&=q\delta_2;&
b_3&=\delta_3;&
b_4&=q^{-1}\delta_3;\\
b_5&=\varepsilon;&
b_6&=q^{-2N-3}\varepsilon;&
g_1&=q^{-2N-4}\varepsilon;&
g_2&=q^{-2N-2}\varepsilon;\\
f_1&=\zeta;&
f_2&=q^{-2N-3}\zeta;&
f_3&=q^{-2N-4}\zeta;&
f_4&=q^{-2N-2}\zeta;\\
&&f_5&=q^{-2N-2}\delta_1;&
a_9&=\delta_1.&&
\end{aligned}
\end{equation}
Moreover, from~\eqref{kl1} we obtain the conditions
\begin{align*}
f_1+a_1+q^2\s+f_3+q^{2N+1}f_2&=0;\\
b_1+a_5+qa_2+q^2\s+e_1+q^{2N+1}a_7&=0;\\
b_3+a_6+a_3+q^2\s+e_2+q^{2N+1}a_8&=0;\\
c+b_5+b_4+qb_2+q^2\s+g_1+q^{2N+1}b_6&=0;\\
g_2+qe_4+e_3+a_9+q^2\s+f_5&=0;
\end{align*}
which lead to further dependencies between the parameters, namely
\begin{equation}\label{p4}
\begin{aligned}
\alpha&=1+(1+q^{-2}+q^{-2N-2}\s+)\zeta;\\
\delta_2&=-\alpha-(1+q^{-2}+q^{-2N-2}\s+)\beta;\\
\delta_3&=-q\alpha-(1+q^{-2}+q^{-2N-2}\s+)\gamma;\\
\varepsilon&=-\beta-q\gamma-(q^2+q^4\s+)\delta_1;\\
c&=-q^2\delta_2-q^{-1}\delta_3-(1+q^{-2}+q^{-2N-2}\s+)\varepsilon,
\end{aligned}
\end{equation}
which we shall leave aside for the moment.
Instead, we evaluate the condition
\begin{equation}\label{kl6}
q\RCPcm stuvklmn \d\x ij\cdot\x st\x uv-\d\x ij\cdot\x kl\x mn=0
\end{equation}
with the system of coefficients reduced only by~\eqref{p9}.  
This leads to
\begin{align*}
\varepsilon&=0;&\beta&=0;&\gamma&=0;&\zeta&=0;&\delta_1&=0,
\end{align*}
thus simplifying \eqref{p4} to
\begin{align*}
\alpha&=1;&\delta_2&=-1;&\delta_3&=-q;&c&=q^2+1.
\end{align*}
Backtracking the substitutions gives the coefficients of the original ansatz
to be
\begin{gather*}
\begin{aligned}
a_1&=0;\qquad&a_2&=q^{-1};\qquad&a_3&=q;\qquad&a_4&=1;\\
b_1&=-1;\qquad&b_2&=-q;\qquad&b_3&=-q;\qquad&b_4&=-1;
\end{aligned}\\
c=q^2+1;\\
a_5=a_6=a_7=a_8=a_9=0;\qquad b_5=b_6=g_1=g_2=0;\\
e_1=e_2=e_3=e_4=0;\qquad f_1=f_2=f_3=f_4=f_5=0.
\end{gather*}
By testing the complete list of necessary conditions
(again with computer-algebra reduction)
one checks that all of them are fulfilled; thus the coefficient system 
obtained describes in fact a differential calculus, with the bimodule structure
being as described in Theorem~\ref{tClass}, case~\ref{it-free}. 
This check is required only in this case since the corresponding statement
for the other two cases is covered by the assertions of Theorem~\ref{tRest}.
This completes the proof of Theorem~\ref{tClass}.\qed

\section{Conclusion}
The description of first order differential calculus on
the quantum projective spaces forms the first step on the way to an
investigation of their noncommutative geometry. Higher order differential
calculus would be the next indispensable pre-requisite for formulating
basic concepts of differential geometry on $\CPq{N-1}$.

At the same time, the results proved here together with previous results
\cite{sch}, \cite{wk1} draw an outline of how covariant first order
differential calculus
on the quantum group $\SUq(N)$ and the related quantum homogeneous spaces
$\Sq^{2N-1}$, $\CPq{N-1}$ is linked and how different closely related
quantum spaces may behave. To illustrate the latter, remember just
that on $\CPq{N-1}$, differential calculi are essentially uniquely
determined while on $\Sq^{2N-1}$ there was a vast variety of parametrical
series of them.

\subsection*{Acknowledgement}
The author gratefully acknowledges support of the work presented here
by Deutsche Forschungsgemeinschaft within the programme of Graduiertenkolleg
Quantenfeldtheorie, Leipzig.

\end{document}